\documentclass[a4paper,11pt]{article}
\usepackage{amsmath}
\usepackage{amssymb}
\numberwithin{equation}{section}
\newtheorem{theo}{Theorem}[section]

\newtheorem{prop}{Proposition}[section]

\newcommand{\Proof}[1]{\noindent{\textit{Proof.}}#1 $\quad\Box$
{\vspace{.2cm}}}
\newcommand{\Rem}[2]{\noindent{\bf Remark #1} {#2} {\vspace{.2cm}}}

\newcommand{\ov}[1]{\overline{#1}}
\newcommand{\pa}{\partial}
\newcommand{\lam}{\lambda}
\newcommand{\ta}{\theta}
\newcommand{\ei}[1]{\mathrm{e}^{#1}}
\newcommand{\C}{\mathbb{C}}
\newcommand{\R}{\mathbb{R}}

\newcommand{\Z}{\mathbb{Z}}

\newcommand{\dis}{\displaystyle}
\newcommand{\om}{\Omega}
\newcommand{\ep}{\epsilon}
\title{Nonrigidity of a class of two dimensional surfaces\\
with positive curvature and planar points}

\author{ Abdelhamid Meziani\\
Department of Mathematics\\
Florida International University\\
Miami, Florida 33199\\
\\
email: meziani@fiu.edu }

\date{}

\begin{document}

\maketitle

\section*{Introduction}
The problem considered here deals with the bendings of an
orientable, embedded surface $S$ in $\R^3$. We assume that $S$ has a
vanishing first homology group, that $\ov{S}$ is a $C^\infty$
compact surface with boundary, that it has positive curvature except
at finitely many planar points in $S$. The main result states that
for any $k\in\Z^+$, $S$ has nontrivial infinitesimal bendings of
class $C^k$.
 That is, there is a $C^k$ function $U:\, \ov{S}\ \longrightarrow\
 \R^3$ such that the first fundamental form of the deformation
 surface
 $S_\sigma =\{ p+\sigma U(p),\ p\in S\}$ satisfies
 $dS_\sigma^2=dS^2+O(\sigma^2)$ as $\sigma\to 0$, where $\sigma$ is
 a real parameter. Furthermore, $S_\sigma$ is not obtained from $S$
 through a rigid motion of $\R^3$. A consequence of this result is
 the nonrigidity of $S$ in the following sense. Any given
 $\ep$-neighborhood of $S$ (for the $C^k$ topology) contains
 isometric surfaces that are not congruent.

 The study of bendings of surfaces in $\R^3$ has a rich history and
 many physical applications. In particular, it is used in the theory
 of elastic shells. We refer to the survey article of
 Sabitov ({\cite{SABITOV1}}) and the references therein.
 The results of this paper are also related to those contained in
 the following papers {\cite{BLEE}}, {\cite{EFIMOV}}, {\cite{GREENE}},
 {\cite{KANN}}, {\cite{KARATO}}, {\cite{MEZCAG}},
 {\cite{MEZJDE}}, {\cite{POGO1}}, {\cite{USMANOV1}}, {\cite{SABITOV2}}

 Our approach is through the study of the associated (complex) field
 of asymptotic directions on $S$. We prove that such a vector field
 generates an integrable structure on $\ov{S}$. We reformulate the
 equations for the bending field $U$ in terms of a Bers-Vekua type
 equation (with singularities). Then use recent results about the
 solvability of such equations to construct the bending fields.

\section{Integrability of the field of asymptotic directions}
For the surfaces considered here, we show that the field of
asymptotic directions on $S$ has a global first integral.

Let $S\subset\R^3$ be an orientable $C^\infty$ surface with a
$C^\infty$ boundary. We assume that $H_1(S)=0$. The surface $S$ is
diffeomorphic to a relatively compact domain $\om\in\R^2$ with a
$C^\infty$ boundary. Hence,
\begin{equation}
\ov{S} =\{ R(s,t)\in\R^3;\ (s,t)\in\ov{\om}\}\, ,
\end{equation}
where the position vector $R:\ \ov{\om}\, \longrightarrow\, \R^3$ is
a $C^\infty$ parametrization of $\ov{S}$. Let $E,\ F,\ G$ and $e,\
f,\ g$ be the coefficients of the first and second fundamental forms
of $S$. Thus,
\[\begin{array}{lll}
E =R_s\cdot R_s,\quad & F =R_s\cdot R_t,\quad & G =R_t\cdot R_t,\\
e =R_{ss}\cdot N,\quad & f =R_{st}\cdot N,\quad & g =R_{tt}\cdot N,
\end{array}\]
where $N=\dis\frac{R_s\times R_t}{|R_s\times R_t|}$ is the unit
normal of $S$. The Gaussian curvature of $S$ is
$K=\dis\frac{eg-f^2}{EG-F^2}$. We assume that $\ov{S}$ has positive
curvature except at a finite number of planar points in $S$. That
is, there exist $p_1=(s_1,t_1)\in\om,\cdots ,p_l(s_l,t_l)\in\om$
such that
\begin{equation}
K(s,t) >0,\quad \forall (s,t)\in\ov{\om}\backslash\{ p_1,\cdots
,p_l\} \, .
\end{equation}
The (complex) asymptotic directions on $S$ are given by the
quadratic equation
\begin{equation}
\lam^2 +2f\lam +eg =0
\end{equation}
Thus $\lam =-f+i\sqrt{eg-f^2}\, \in\, \R+i\R^+$ except at the planar
points $p_1,\cdots ,p_l$ where $\lam =0$.

Consider the structure on $\ov{\om}$ generated by the $\C$-valued
vector field
\begin{equation}
L=g(s,t)\frac{\pa}{\pa s}+\lam(s,t)\frac{\pa}{\pa t}\, .
\end{equation}
This structure is elliptic on $\ov{\om}\backslash\{ p_1,\cdots
,p_l\} $. That is, $L$ and $\ov{L}$ are independent outside the
planar points. The next proposition shows that $L$ has a global
first integral on $\ov{S}$.

{\prop{ Let S be surface given by $(1.1)$ whose curvature $K$
satisfies $(1.2)$. Then there exists an injective function
 \[
Z:\ \ov{\om}\, \longrightarrow\, \C
 \]
such that
\begin{enumerate}
\item[$1.$] $Z$ is $C^\infty$ on $\ov{\om}\backslash\{ p_1,\cdots
,p_l\} $;
\item[$2.$] $LZ=0$ on $\ov{\om}\backslash\{ p_1,\cdots
,p_l\} $; and
\item[$3.$] For every $j=1,\cdots ,l$, there exists $\mu_j >0$ and polar
coordinates $(r,\ta)$ centered at $p_j$ such that in neighborhood of
$p_j$ we have
\begin{equation}
Z(r,\ta)=Z(0,0)+r^{\mu_j}\ei{i\ta}+O(r^{2\mu_j})
\end{equation}
\end{enumerate}
 }}

 \Proof{ Since L is $C^\infty$ and elliptic on
$\ov{\om}\backslash\{ p_1,\cdots ,p_l\} $, then it follows from the
uniformization of complex structures on planar domains (see{\cite{
SPRINGER}}) that there exists a $C^\infty$ diffeomorphism
\[
Z: \ \ov{\om}\backslash\{ p_1,\cdots ,p_l\}\,\longrightarrow\,
Z(\ov{\om}\backslash\{ p_1,\cdots ,p_l\} )\subset\C
\]
such that $LZ=0$. It remains to show that $Z$ has the form (1.5) in
a neighborhood of a planar point.

Let $p_j$ be a planar point of $S$. We can assume that $S$ is given
in a neighborhood of $p_j$ as the graph of a function $z=z(x,y)$
with $p_j=(0,0)$, $z(0,0)=0$, and $z_x(0,0)=z_y(0,0)=0$. The
assumption on the curvature implies that
 \[
 z(x,y)=z_m(x,y) +o(\sqrt{x^2+y^2}^m),
 \]
 where $z_m(x,y)$ is a homogeneous polynomial of degree $m>2$,
 satisfying $z_{xx}z_{yy}-z_{xy}^2 >0$ for $(x,y)\ne 0$. We can also
 assume that $z(x,y) >0$ for $(x,y)\ne 0$. The complex structure
 generated by the asymptotic directions is given by the vector field
 \[
 L=z_{yy}\frac{\pa}{\pa x}+(-z_{xy}+i\sqrt{z_{xx}z_{yy}-z_{xy}^2})
 \frac{\pa}{\pa y}\, .
 \]
 With respect to the polar coordinates $x=\rho\cos\phi$,
 $y=\rho\sin\phi$, we get
 \[
 z=\rho^mP(\phi)+\rho^{m+1}A(\rho,\phi)\, ,
 \]
where $P(\phi)$ is a trigonometric polynomial of degree $m$
satisfying $P(\phi)>0$ and (curvature)
\[
 m^2P(\phi)^2+mP(\phi)P''(\phi)-(m-1)P'(\phi)^2 > 0\,\qquad
 \forall\phi \in\R\, .
\]
With respect to the coordinates $(\rho,\phi)$, the vector field $L$
becomes
\[
L=m(m-1)\rho^{m-2}(P(\phi)+O(\rho))L_0\,
\]
with
\[
L_0=\frac{\pa}{\pa\phi}+\rho
\left(M(\phi)+iN(\phi)+O(\rho)\right)\frac{\pa}{\pa\rho}
\]
and
 \[
 M=\frac{P'}{mP}\quad\mathrm{and}\quad
 N=\frac{1}{m}\sqrt{\frac{m^2P^2+mPP''-(m-1)P'^2}{(m-1)P^2}}\, .
 \]
 We know (see{\cite{MEZJFA}}) that such a vector field $L_0$ is
 integrable in a neighborhood of the circle $\rho =0$. Moreover, we
 can find coordinates $(r,\ta)$ in which $L_0$ is $C^1$-conjugate to
 the model vector field
 \[
T=\mu_j \frac{\pa}{\pa\ta}-ir\frac{\pa}{\pa r}
 \]
where $\mu_j >0$ is given by
\[
\frac{1}{\mu_j}=\frac{1}{2\pi}\int_0^{2\pi}(N(\phi)-iM(\phi))d\phi
=\frac{1}{2\pi}\int_0^{2\pi}N(\phi)d\phi\, .
\]
The function $u_j(r,\ta)=r^{\mu_j}\ei{i\ta}$  is a first integral of
$T$ in $r>0$.

Now we prove that the function $Z$ which is defined in
$\ov{\om}\backslash\{ p_1,\cdots ,p_l\}$ extends to $p_j$ with the
desired form given by (1.5). Let $O_j$ be a disc centered at $p_j$
where $L$ is conjugate to a multiple of $T$ in the $(r,\ta)$
coordinates. Since $u_j$ and $Z$ are both first integrals of $L$ in
the punctured disc $O_j\backslash p_j$, then there exists a
holomorphic function $h_j$ defined on the image $u_j(O_j\backslash
p_j)$ such that
 $ Z(r,t)=h_j(u_j(r,t))$. Since both $Z$ and $u_j$ are
 homeomorphisms onto their images, then $h_j$ is one to one in a
 neighborhood of $u_j(p_j)=0\in\C$ and since $h_j$ is bounded, then
 \[
 h_j(\zeta )=C_0+C_1\zeta +O(\zeta^2)\qquad\mathrm{for}\ \zeta
 \ \mathrm{close\ to\ }0\in\C
 \]
with $C_1\ne 0$. This means that after a linear change of the
coordinates $(r,\ta)$ (to remove the constant $C_1$) , the function
$Z$ has the form (1.5) }

\section{Equations of the bending fields in terms of $L$}
Let $S$ be a surface given by (1.1). An infinitesimal bending of
class $C^k$ of $S$ is a deformation surface
$S_\sigma\subset\in\R^3$, with $\sigma\in\R$ a parameter, given by
the position vector
\begin{equation}
R_\sigma (s,t)=R(s,t)+\sigma U(s,t)\, ,
\end{equation}
whose first fundamental form satisfies
\[
 dR_\sigma^2=dR^2+O(\sigma^2)\qquad\mathrm{as}\ \sigma\to 0\, .
\]
This means that the bending field $U:\, \ov{\om}\, \longrightarrow\,
\R^3$ is of class $C^k$ and satisfies
\begin{equation}
dR\cdot dU =0\, .
\end{equation}
The trivial bendings of $S$ are those induced by the rigid motions
of $\R^3$. They are given by $U(s,t)=A\times R(s,t)+B$, where $A$
and $B$ are constants in $\R^3$, and where $\times$ denotes the
vector product in $\R^3$.

Let $L$ be the field of asymptotic directions defined by (1.4). For
each function $U:\, \om\,\longrightarrow\, \R^3$, we associate the
$\C$-valued function $w$ defined by
\begin{equation}
w(s,t)=LR(s,t)\cdot U(s,t)=g(s,t)u(s,t)+\lam (s,t)v(s,t)\, ,
\end{equation}
where $u=R_s\cdot U$ and $v=R_t\cdot U$. The following theorem
proved in {\cite{MEZCONT}} will be used in the next section.

{\theo{\rm{\cite{MEZCONT}}} If $U:\, \om\,\longrightarrow\,\R^3$
satisfies $(2.2)$, then the function $w$ given by $(2.3)$ satisfies
the equation
 \begin{equation}
 CLw=Aw+B\ov{w},
 \end{equation}
where
\begin{equation}\begin{array}{ll}
 A & = (LR\times\ov{L}R)\cdot (L^2R\times\ov{L}R)\, ,\\
 B & = (LR\times\ov{L}R)\cdot (L^2R\times LR)\, ,\\
 C & = (LR\times\ov{L}R)\cdot (LR\times\ov{L}R)\, .
\end{array}\end{equation}
}

\Rem{2.1}{ If $w$ solves equation (2.4). The function $w'=aw$, where
$a$ is a nonvanishing function solves the same equation with the
vector field $L$ replaced by $L'=aL$}

\section{Main Results}

{\theo Let $S$ be a surface given by $(1.1)$ and such that its
curvature $K$ satisfies $(1.2)$. Then for every $k\in\Z^+$, the
surface $S$ has a nontrivial infinitesimal bending $U:\,
\ov{\om}\,\longrightarrow\, \R^3$ of class $C^k$.}

\vspace{.2cm}

 \Rem{3.1}{ It should be mentioned that without the
assumption that $K>0$ up to the boundary $\pa S$, the surface could
be rigid under infinitesimal bendings. Indeed, let $T^2$ be a
standard torus in $\R^3$, it is known (see{\cite{AUDOLY}} or
{\cite{POGO2}}) that if $S$ consists of the portion of $T^2$ with
positive curvature, then $S$ is rigid under infinitesimal bendings.
Here the curvature vanishes on $\pa S$.}

Before we proceed with the proof, we give a consequence of Theorem
3.1.

{\theo Let $S$ be as in Theorem 3.1. Then for every $\ep >0$ and for
every $k\in\Z^+$, there exist surfaces $\Sigma^+$ and $\Sigma^-$ of
class $C^k$ in the $\ep$-neighborhood of $S$ (for the
$C^k$-topology) such that $\Sigma^+$ and $\Sigma^-$ are isometric
but not congruent.}

\vspace{.2cm}

\Proof{ Let $U:\, \ov{\om}\,\longrightarrow\, \R^3$ be a nontrivial
infinitesimal bending of $S$ of class $C^k$. Consider the surfaces
$\Sigma_{\sigma}$ and $\Sigma_{-\sigma}$ defined the position
vectors
\[
R_{\pm\sigma}(s,t)=R(s,t)\pm\sigma U(s,t)\, .
\]
Since $dR\cdot dU=0$, then $dR_{\pm\sigma}^2=dR^2+\sigma^2dU^2$.
Hence $\Sigma_\sigma$ and $\Sigma_{-\sigma}$ are isometric.
Furthermore, since $U$ is nontrivial, then $\Sigma_\sigma$ and
$\Sigma_{-\sigma}$ are not congruent (see {\cite{SPIVAK}}). For a
given $\ep >0$, the surfaces $\Sigma_{\pm\sigma}$ are contained in
the $\ep$-neighborhood of $S$ if $\sigma$ is small enough }

\noindent{\it Proof of Theorem 3.1.} First we construct non trivial
solutions $w$ of equation (2.4) and then deduce the infinitesimal
bending fields $U$. For this, we use the first integral $Z$ of $L$
to transform equation (2.4) into a Bers-Vekua type equation with
singularities. Let $Z_1=Z(p_1),\cdots , Z_l=Z(p_l)$ be the images of
the planar points by the function $Z$. The pushforward of equation
(2.4) via $Z$ gives rise to an equation of the form
\begin{equation}
\frac{\pa W}{\pa\ov{Z}}=\frac{A(Z)}{\prod_{j=1}^l(Z-Z_j) }W+
\frac{B(Z)}{\prod_{j=1}^l(Z-Z_j) }\ov{W}\, ,
\end{equation}
where $w(s,t)=W(Z(s,t))$ and
$
 A,B\, \in\, C^\infty(Z(\ov{\om}\backslash\{ p_1,\cdots
 ,p_l\}))\cap L^\infty (Z(\ov{\om})).
$
 The local study of the solutions of such equations near a
 singularity is considered in {\cite{MEZCVEE}}, {\cite{TUNG}},
 {\cite{USMANOV2}}. To construct a global solution of (3.1)
 with the desired properties, we proceed as follows. We seek a
 solution $W$ in the form $W(Z)=H(Z)W_1(Z)$, where
 \[
 H(Z)=\prod_{j=1}^l(Z-Z_j)^M\, ,
 \]
where $M$ is a (large) positive integer to be chosen. In order for
$W=HW_1$ to solve (3.1), the function $W_1$ needs to solve the
modified equation
\begin{equation}
\frac{\pa W_1}{\pa\ov{Z}}=\frac{A(Z)}{\prod_{j=1}^l(Z-Z_j) }W_1+
\frac{B(Z)}{\prod_{j=1}^l(Z-Z_j) }\frac{\ov{H(Z)}}{H(Z)}\ov{W_1}\, .
\end{equation}
Since $A$ and $B\ov{H}/H$ are bounded functions on $Z(\ov{\om})$,
then a result of {\cite{TUNG}} gives a continuous solution $W_1$ of
(3.2) on $Z(\ov{\om})$.  Furthermore, such a solutions is $C^\infty$
on $Z(\ov{\om}\backslash\{Z_1,\cdots ,Z_l\})$ since the equation is
elliptic  and the coefficients are $C^\infty$ outside the $Z_j$'s.

The function $W(Z)=H(Z)W_1(Z)$ is therefore a solution of (3.1). It
is $C^\infty$ on $Z(\ov{\om}\backslash\{Z_1,\cdots ,Z_l\})$ and
vanishes to order $M$ at each point $Z_j$. Consequently, the
function $w(s,t)=W(Z(s,t))$ is $C^\infty$ on
$\ov{\om}\backslash\{p_1,\cdots ,p_l\}$ and vanishes to order
$M\mu_j$ at each planar point ($\mu_j$ is the positive number
appearing in Proposition 1.1).

Now we recover the bending field $U$ from the solution $w$ of (2.4)
and the relation $w=LR\cdot U$. Set $w=gu+\lam v$, where $\lam$ is
the asymptotic direction given in (1.3). The functions $u$ and $v$
are uniquely determined by
\[
 v=\frac{w-\ov{w}}{2i\sqrt{eg-f^2}}\quad\mathrm{and}\quad
 u=\frac{w+\ov{w}+2fv}{2g}\, ,
\]
provided that the function $w$ vanishes to a high order at the
planar points (order of vanishing of $W$ at $p_j$ larger than that
of the curvature). These functions are $C^\infty$ outside the planar
points. At each planar point $p_j$, it $m_j$ is the order of
vanishing of $K$, then the functions $v$ and $u$ vanish to order
$M\mu_j-m_j$. It follows from $LR\cdot U=w$ that $R_s\cdot U=u$ and
$R_t\cdot U=v$. The condition $dR\cdot dU=0$ implies that
 \begin{equation}
R_{ss}\cdot U=u_s,\quad R_{tt}\cdot U=v_t,\quad\mathrm{and}\quad
2R_{st}\cdot U=u_t+v_s\, .
 \end{equation}
In terms of the components $(x,y,z)$ of $R$ and $(\xi,\eta,\zeta)$
of $U$, we have
\begin{equation}
\left\{\begin{array}{ll}
 x_s\xi+y_s\eta+z_s\zeta & =u\\
 x_t\xi +y_t\eta +z_t\zeta & =v\\
 x_{ss}\xi +y_{ss}\eta+z_{ss}\zeta & =u_s\\
 x_{tt}\xi+y_{tt}\eta +z_{tt}\zeta & =v_t\\
 2x_{st}\xi+2y_{st}\eta+2z_{st}\zeta & =u_t+v_s
\end{array}\right.
\end{equation}
(equation (2.4) guarantees the compatibility of this system). Note
that at each point $p\in\ov{\om}$ where $K>0$, the functions $\xi,\
\eta,$ and $\zeta$ are uniquely determined by $u$, $v$, and $u_s$
(or $v_t$). Indeed, at such a point the determinant of the first
three equations of (3.4) is
\[
 R_{ss}\cdot (R_s\times R_t)=|R_s\times R_t|e\, \ne \, 0\, .
\]
With our choice that $w$ (and so $u$ and $v$) vanishing to an order
larger than that of the curvature at each planar point, the
functions $\xi,\ \eta,$ and $\zeta$ are also uniquely determined to
be 0 at each planar point. To see why, assume that at $p_j$, we have
$x_sy_t-x_ty_s\ne 0$, then after solving the first two equations for
$\xi$ and $\eta$ in terms of $\zeta$, $u$, and $v$, the third
equations becomes,
\begin{equation}
 R_{ss}\cdot (R_s\times R_t) \zeta =\left|
\begin{array}{lll}
 x_s & y_s & u\\
 x_t & y_t & v\\
 x_{ss} & y_{ss} & u_s
\end{array}
 \right| .
\end{equation}
Since the zeros $p_j$ of $R_{ss}\cdot (R_s\times R_t)$ is isolated
and since $u$ and $v$ vanish to a high order at $p_j$, the function
$\zeta$ is well defined by  (3.5). Consequently, for any given
$k\in\Z^+$, a nonzero solution $w$ of (2.4) which vanishes at high
orders ($M$ large), gives rise to a unique field of infinitesimal
bending $U$ of $S$,  so that it is $C^\infty$ on
$\ov{\om}\backslash\{p_1,\cdots ,p_l\}$ and vanishes to an order $k$
at each $p_j$. Such a field is therefore of class $C^k$ at each
planar point. It remains to verify that $U$ is not trivial. If such
a field were trivial ($U=A\times R +B$), then the vanishing of
$dU=A\times dR$ at $p_j$ together with $dR\cdot dU=0$ gives $A=0$
and so $U=B=0$ since $U=0$ at $p_j$. This would give $w\equiv 0$
which is a contradiction $\quad\Box$

\end{document}